\documentclass[10pt]{article}
\usepackage{graphicx}
\usepackage{amssymb}
\usepackage{psfrag}
\usepackage{soul}
\usepackage{multirow}
\usepackage{makecell}
\title{Tessellations}
\author{Chaim Goodman-Strauss}
\date{}

\def\ddefine#1{{\em #1}}

\def\qed{\hfill{$\Box$}}

\def\vp{\vspace{.5\baselineskip}}

\newenvironment{proof}{\vp 
   \par \noindent {\bf Proof}\ } {\hfill $\Box$ \medskip \par}

\def\compformat{\tt}

\def\mp#1{\marginpar{\tiny\em #1}}

\def\vect#1{{\bf #1}}

\begin{document}

\maketitle

{\em \small (This article appeared in Italian as {\em Tassellazioni}, in La Matematica {\bf 3}, Einaudi (2011), Turin, pages 249-285.})

\vp

\begin{figure}\centerline{\includegraphics[width=.8 \textwidth]{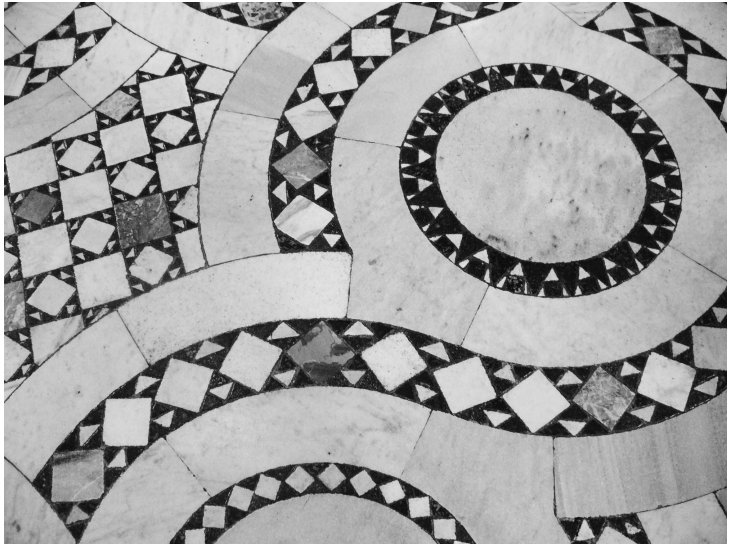}}
\caption{The late cosmatesque floors of San Giovanni in Laterano,  completed in 1425.}\label{floor}\end{figure}

Tessellations --- the patterned repetition of small ceramic or stone pieces ---  appear all over the world, in virtually every culture with a tradition of permanent, formal building, over the last several thousand years, perhaps reaching an apotheosis in the splendid and sumptuous traditions of the Islamic world. The technique splendidly solves an important design question: how to create a pleasing ornamental surface of any needed size, as in Figure~\ref{floor}. 

The craft of tile work, {\it per se}, requires formally constructed buildings with permanent surfaces that can be tessellated, structures only possible in more complex civilizations.  However, human preoccupation with geometric ornament is far older still, and far more widespread. In virtually every human society, textiles, wall ornament and pottery illuminate a fundamental human desire for the pleasing repetition of motifs: one can ask for no better study than Dorothy Washburn and Donald Crowe's {\em Symmetries of Cultures}~\cite{washcrowe}.

Truly archaic examples are difficult to identify, perhaps because materials such as fiber, wood or leather are preserved only in the most unusual circumstances. But the remarkable
 ochres found in the Blombos caves, carved some 77,000 years ago, show unmistakable evidence of geometrical thought. Exploring pattern and symmetry are truly as human an endeavor as can be, as central to the human experience as language and rhythm.
 
 This essay is appearing in a volume on the ties of mathematics to the broader culture, and exploring this theme in some depth would be natural, perhaps expected. But I would like to take a somewhat different approach, and discuss something of the {\em contemporary} and {\em mathematical} culture of tessellations.  
 
 Perhaps because of the subject's deep roots in the human experience, the study of tessellations links from the arts through recreational mathematics all the way to one of the deepest threads in contemporary mathematics, the theory of computation, the question of what can and what cannot be computed.  It is notable that many of the remarkable constructions in this essay are not due to professional academics, but to, well, there is no better description than {\em enthusiasts!} Few areas of contemporary mathematics are so amenable to important contributions from the non-professional--- we've included a number of open questions that perhaps {\em you} will answer for the rest of us! 
 
 If you are interested in exploring further, you will find no more comprehensive source than Branko Gr\"unbaum and Geoffrey Shepherd's masterpiece, {\em Tilings and Patterns}~\cite{grsh}. Doris Schattschneider's {\em M.C. Escher: Visions of Symmetry}~\cite{schattschneider} illuminates the tessellation work of the inimitable artist, who anticipated many themes later explored by academic mathematicians. {\em The Symmetries of Things} demonstrates contemporary topological tools for examining these issues~\cite{sot}. The many works of Martin Gardner, and the volume of essays, {\em Mathematical Recreations}~\cite{klarner}, in his honor, are very useful to the young tessellator. Finally, most importantly, overshadowing all of these, are the countless works by innumerable artists and artisans through time and space, from whom this mathematical inquiry draws inspiration.
 
 In this essay, we won't try to survey every possible interesting corner--- there's just too much!--- but set our to explore an ancient question, somewhat formalized here, to be sure: {\em What can fit together, and how?} and follow a particular thread that  leads us to one of the deeper issues in contemporary mathematics, that some problems are {\em undecidable}, can be {\em proven} to be unanswerable by mechanical means.\marginpar{\tiny cross ref elsewhere in La Matamatica}



%


%


\psfrag{A}{$A$}
\psfrag{B}{$B$}
\psfrag{C}{$C$}
\psfrag{D}{$D$}
\psfrag{E}{$E$}

\section{What can tessellate?}\label{liviosec}

\begin{figure}[h]\centerline{\includegraphics{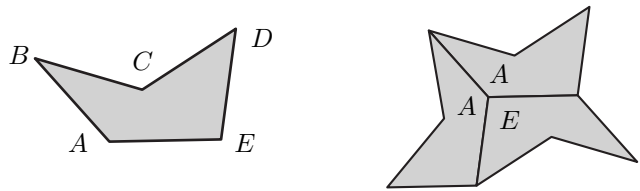}}\caption{Livio's pentagon.}\label{livio0}\end{figure}

Let us start with the non-convex, equilateral pentagon, discovered by Livio Zucci (Figure~\ref{livio0}),  the unique equilateral pentagon for which the vertex angles satisfy $2A+E= B+C+2D=B+C+E=2\pi$, describing some of the ways copies of the tile can fit together around a vertex (for example, two copies of the $A$ vertex can fit with an $E$ vertex, as shown). These relationships imply several others, and copies of the tile can form a remarkable variety of configurations and patterns.  You really must have a set of these tiles to play with in order to appreciate their rich behavior  ---  I suggest buying several copies of this volume, saving one, and cutting out sets of tiles from all the rest.\footnote{Alternatively, materials can be downloaded from the author's website {\em http://mathfactor.uark.edu/downloads/tiles.pdf} 
}

\begin{figure}\centerline{\includegraphics{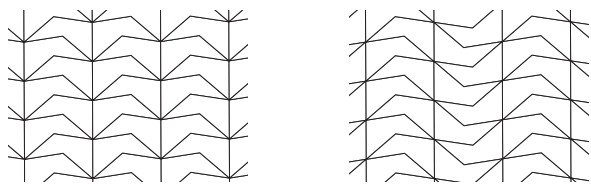}}\caption{Two ``isohedral" tessellations admitted by Livio's pentagon.}\label{livio5}\end{figure}

It's quite a subtle problem: just what tessellations {\em can} this tile form-- and what {\em can't} it? In Figure~\ref{livio5} are two patterns that can each be extended to \ddefine{isohedral} tessellations ---  tessellations with a symmetry that acts transitively on the tiles themselves. In more familiar language, in such a tessellation, each tile has the same relationship to the whole as every other tile does.

\begin{figure}
\centerline{\includegraphics[width = 4cm]{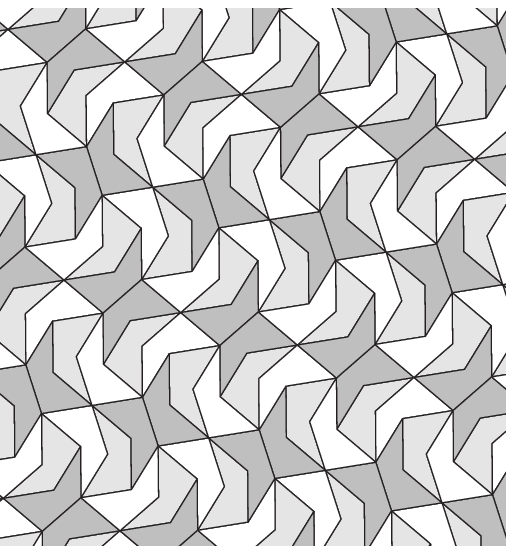}\hspace{\baselineskip}{\includegraphics[width = 4cm]{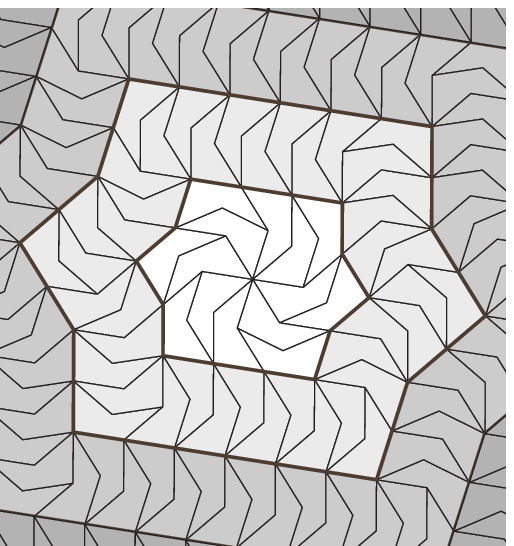}}}

\vspace{\baselineskip}
\centerline{\includegraphics[width = 4cm]{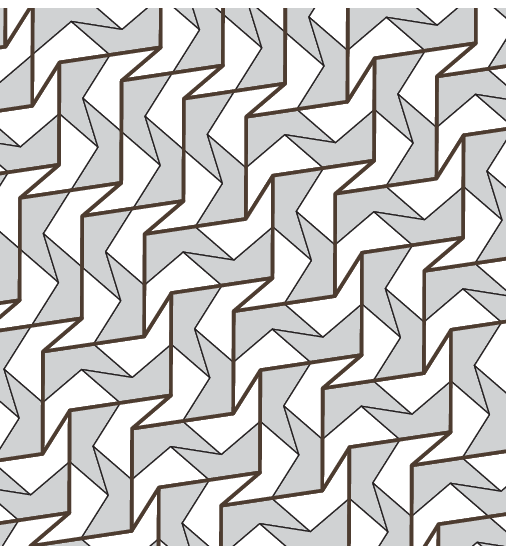}\hspace{\baselineskip}{\includegraphics[width = 4cm]{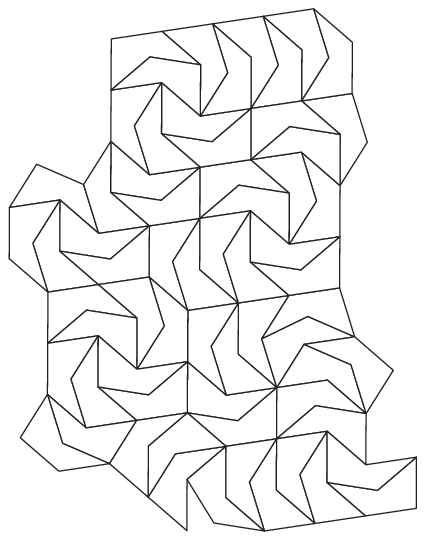}}}
\caption{Three ``anisohedral" tessellations by Livio's pentagon. 
 Can the configuration
 at bottom right be completed into a tessellation of the entire Euclidean plane? More importantly, how might one tell?}
\label{livio}
\end{figure}

In Figure~\ref{livio} we begin to see the rich behavior of Livio's pentagon. 
It is clear that three of the patterns in the figure can be extended to \ddefine{anisohedral} tessellations of the entire plane, the tiles lying in more than one aspect with relation to the whole, the  symmetry  not acting transitively on the tiles.  If we continue the pattern shown at upper left across the infinite Euclidean plane, we will have a  \ddefine{periodic} tessellation ---  that is, one formed of regularly repeating units (or ``fundamental domains"). Another way to put this is that if you were to turn your back for a moment, a friend could slide the tessellation some distance, yet when you look again you would be unable to tell that anything had changed. On the other hand, the tessellation is still anisohedral ---  the tiles lie in a few different aspects relative to the whole tessellation, indicated by the way they are colored. 

Though it exhibits a definite pattern, the tessellation with concentric rings of tiles, shown at upper right, is \ddefine{non-periodic}: Because the rings are all centered on a particular point, if your friend were to shift the tessellation in any way, you would certainly notice. Is the figure at lower left a portion of a periodic or of a non-periodic tessellation? No finite figure can really tell us ---  this pattern can be completed  either way, but if this tessellation is to be periodic, the units must be fairly large, larger than the image shown here.  

But what about  the \ddefine{configuration} ---  or tessellation of a disk-like region ---  at the lower lower right of the figure? Can it be extended to a tessellation of the entire plane? Most important of all, what criteria could we apply to work this out?


%

\section{What tiles admit tessellations?}\label{whattiles}

The question before us, the primary question underlying the study of tessellations, is how do local constraints inform global structure? As Livio's example richly shows, the question is remarkably subtle.

Given a tile, we can ask the most basic question of all: {\em Does the tile admit a tessellation of the plane, that is, does there exist a tessellation of the plane formed of copies of the tile, with no gaps or overlaps?} More to the point: how can we tell? As we shall learn, there are surprising and  deep lacuna in our understanding. 

\begin{figure}
\centerline{\includegraphics[width = 4cm]{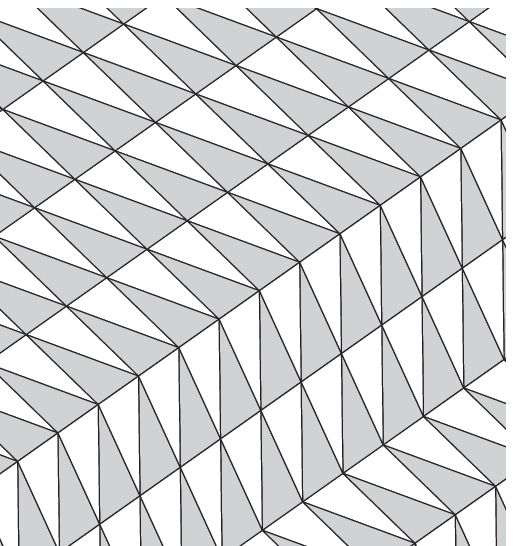}\hspace{\baselineskip}{\includegraphics[width = 4cm]{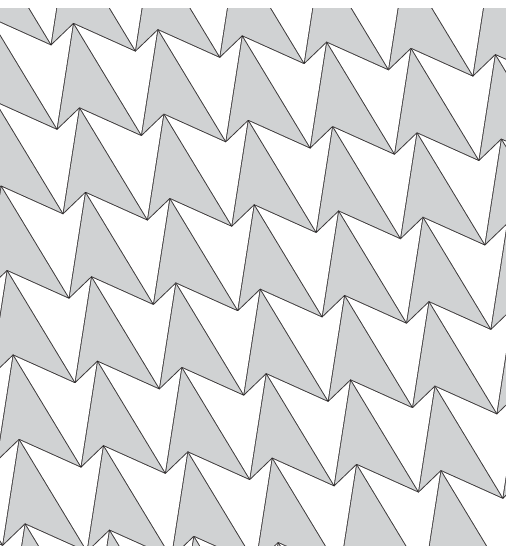}}}
\caption{Each triangle and each quadrilateral admits at least some tessellation.}
\label{triquad}
\end{figure}

But to begin with, as suggested by Figure~\ref{triquad}, it is not difficult to see that every triangle admits a tiling of the plane, as does every quadrilateral, convex or non-convex, though it is considerably more subtle to ask about all tilings admitted by a given triangle or quadrilateral!

The Euclidean Tiling Theorem, which we state and prove at the end of this essay, provides some powerful constraints on the number of neighbors a tile can have in a tessellation of the Euclidean plane. In particular, in a tessellation by copies of a single tile, the Euclidean Tiling Theorem tells us that it is impossible for every copy to have seven or more neighbors. 

Now if the tile is {\em  convex}, then every copy will have at least as many neighbors as it has sides, and so we can go further and be sure that if a convex tile admits a tessellation then it has six or fewer sides.

 \begin{figure}
\centerline{\includegraphics[width = 5cm]{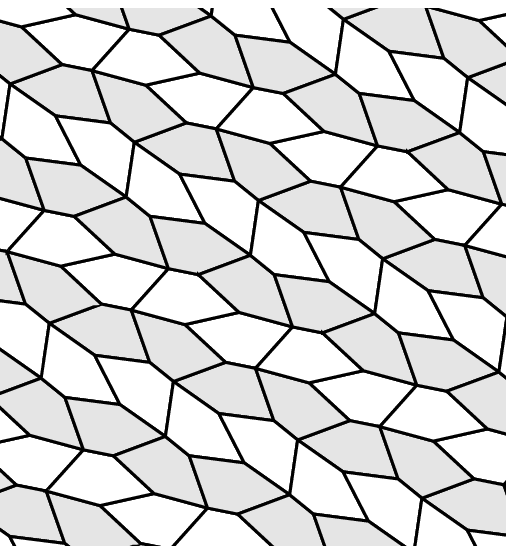}\hspace{\baselineskip}{\includegraphics[width = 5cm]{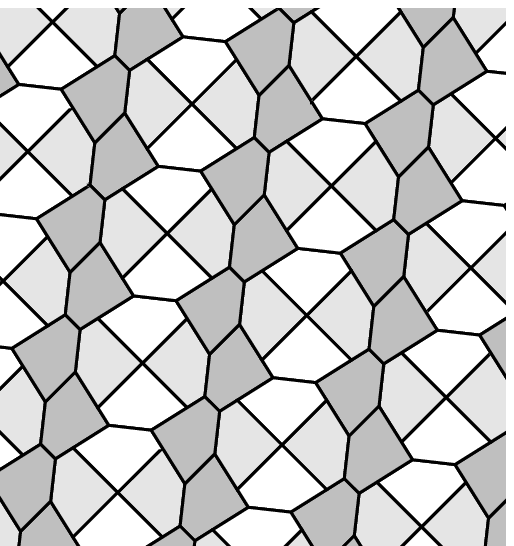}}}

\vspace{\baselineskip}
\centerline{\includegraphics[width = 5cm]{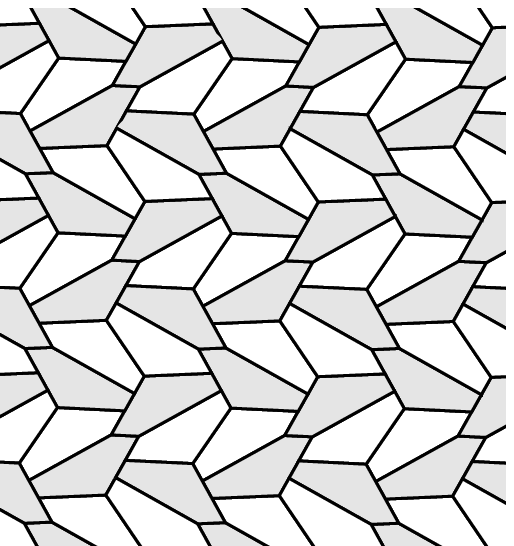}\hspace{\baselineskip}{\includegraphics[width = 5cm]{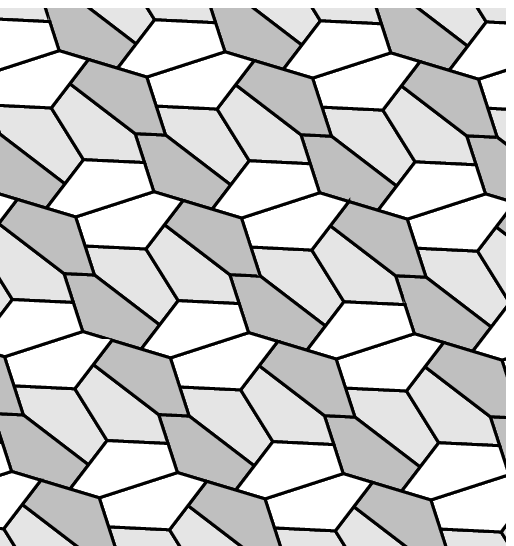}}}
\caption{Four periodic tessellations by anisohedral pentagonal tiles, discovered (clockwise from upper left) by Kersher, Rice, James and Stein, with isohedral numbers two, three, three and two respectively ---  that is, for example, in any tessellation by the pentagon at upper left the tiles are in a minimum of two aspects with respect to the whole. 
Still today, there is very little understanding of which convex pentagons admit tessellations, much less how to address the Monotiling Problem on page [[{\em cross reference}]].}
\label{penta}
\end{figure}
 
We've already seen that every quadrilateral and every triangle admits a tiling, so the remaining questions here are: Which convex pentagons admit tessellations?  Which convex hexagons admit tessellations?\footnote{As with so much of the topic of tessellations, Gr\"unbaum and Shephard's {\em Tilings and Patterns}~\cite{grsh} remains the definitive account. We refer the reader to their Chapter 9 for a more detailed account of the history of these questions and for further bibliographic references, many of which appear in a wonderful collection of articles in honor of Martin Gardner~\cite{klarner}}

If we add one more constraint, that the tessellations be isohedral,  as those in Figure~\ref{livio5} are, these questions are quite easily answered: It is not particularly difficult to enumerate {\em all} isohedral tessellations, as in~\cite{grsh}, or with more topological methods in~\cite{sot}, checking as we go to see which of these are by pentagons or hexagons, whether convex or not. Doing so, we find there are exactly five families of pentagons and three families of hexagons that admit isohedral tessellations. 

The topic is instantly richer when we consider {\em anisohedral} tiles, tiles that admit tessellations, but only anisohedral ones. In 1918, K. Reinhardt showed there can be no convex anisohedral hexagon, but starting in the late 1960's, a number of convex anisohedral  pentagons began to be discovered, by Richard Kershner, then later by Richard James III, Marjorie Rice and finally Rolf Stein.\footnote{Since this essay appeared in {\em La Matematica}, a fifteenth pentagon, with isohedral number 3, was found by Casey Mann, Jennifer McLoud, and David Von Derau.  See {\em 
arXiv:1510.01186}} A few of these are illustrated in Figure~\ref{penta} and an interactive {\em Mathematica} notebook can be obtained at~\cite{pegg}.  But to this day, the question remains quite open: {\em What are the convex anisohedral pentagons?}

We can go further and ask for non-convex anisohedral hexagons and pentagons. Though there has been much investigation and many examples are known, no fully general picture has emerged.

Indeed, it can be remarkably difficult to determine whether a given tile does in fact admit a tessellation of the plane. This clearly depends on the ways copies of the tile can fit together, but even within a small region, this problem is full of surprises. For example, in 1936, H. Voderberg gave a remarkable construction for any given $n$, producing a tile for which two copies can surround $n$ others! His tiles don't admit tessellations for $n>2$, but there is no reason, at present, to rule out the possibility of tiles that have this property, and do tile the plane. Analyzing, then, what can and cannot occur, begins to seem quite subtle. Many other interesting examples can be found throughout~\cite{grsh}. We pose:

 \begin{figure}\centerline{\includegraphics[width=.8\textwidth]{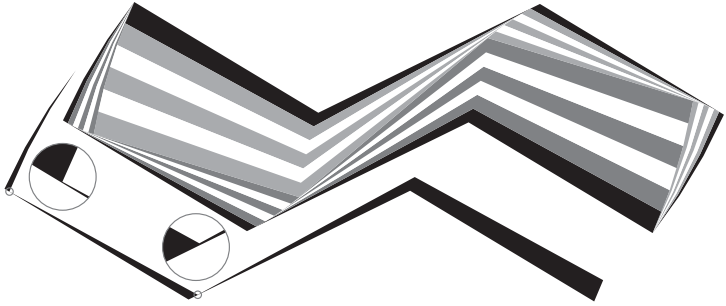}}\caption{Voderberg's tile, here with $n=10$: ten (!) copies are surrounded by two others. Careful inspection reveals the tile really is in just one piece.}\label{voderberg}\end{figure}

\vspace{\baselineskip}\noindent{\bf The Monotiling Problem: } {\em Does a given tile admit a tessellation of the plane?}

\vspace{\baselineskip} 

For any of the tiles we have seen so far, this is a relatively easy problem to settle, even if the tiles have strange properties. With just a little experimentation (or a glance at Figure~\ref{livio}) we can decide that Livio's Pentagon does admit a tessellation, or that the Voderberg tile in Figure~\ref{voderberg} does not. 

\begin{figure}\centerline{\includegraphics[]{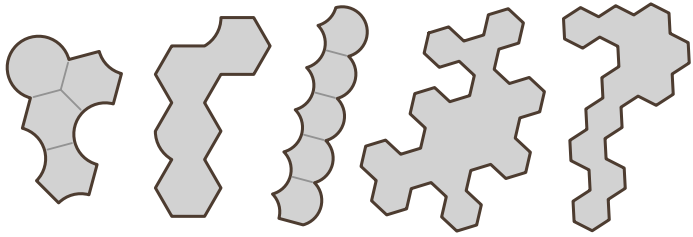}}\caption{The first three tiles are among several discovered by Casey Mann in 2000 and 2007; the two on right are among several discovered by Joseph Myers in 2003. Which of these admits a tessellation of the plane?}\label{mannmyers}\end{figure}

But what about the remarkable tiles, discovered by Casey Mann and Joseph Myers, shown in Figure~\ref{mannmyers}? Which of these does or does not admit a tessellation?

It is not too difficult to prove, at least, that the tile at left does not admit a tessellation: it can be viewed as a cluster of hexagons, with a total of seven edges bulging inwards and four bulging outwards. In any configuration of $N$ tiles, there are at least $3N$ unmatched inward-bulging edges, which must lie on the boundary of the configuration. If the tile were to admit a tessellation, we could cover disks of arbitrary radius $R$;   the number of tiles in such a configuration, and hence the number of unmatched inward-bulging edges, would be approximately proportional to $R^2$, but the number of the edges on the boundary would be approximately proportional to only the circumference of the disk and hence to $R$. For small radii, this may be possible, but asymptotically, at some point, we have too many tiles and not enough room for excess unmatched edges. The tile cannot admit a tessellation. 

We can similarly prove the third tile does not admit a tessellation, but this will not work for the second: it seems we must enumerate all possible configurations, covering  disks of larger and larger radii. Eventually, we run out of configurations, and will have proven that this tile cannot admit a tessellation. 

 \begin{figure}\centerline{\includegraphics[]{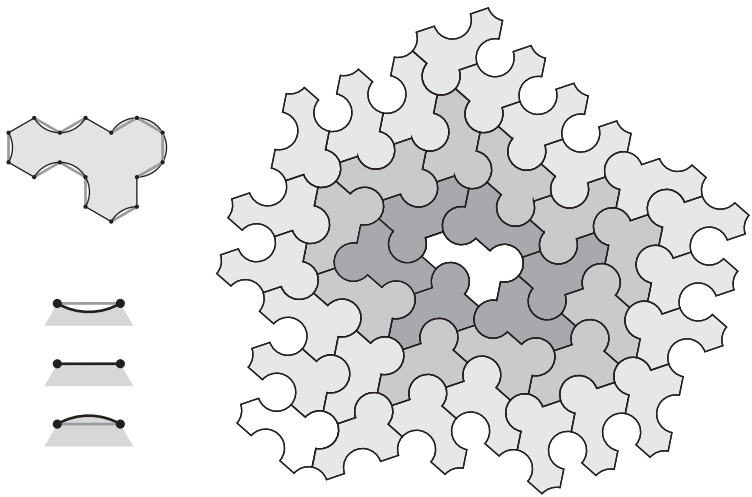}}\caption{This tile does not admit any tessellation of the plane,  and has Heesch number 3 --- a fourth ring can be started, as at lower left of the configuration, but not completed \cite{mann07}.}\label{mann}\end{figure}

The \ddefine{Heesch number} measures the complexity of such a tile, as the largest ``combinatorial radius" of disks it can cover ---  that is the number of concentric rings that copies of the tile can form. In Figure~\ref{mann} we see the first tile has Heesch number 2. The current world record, incidentally, is held by the tile at the center of Figure~\ref{mannmyers}, with Heesch number 5~\cite{mannHeesch}. But still today, no one knows the answer to the question: {\em For all $n$, is there a tile with Heesch number at least $n$?}

The two tiles at right {\em do} admit tessellations, though I suspect that this is humanly impossible to discover! (I'd very much like to know if it is!) But at least there is a mechanical means available:  If a tile admits a {\em periodic} tessellation, we can find this out  using another brute force enumeration of  all possible configurations, this time checking as we go whether or not we have one that can be used as a fundamental domain. The isohedral number, defined in Figure~\ref{penta} above, measures the complexity of this.

\begin{figure}
\centerline{\includegraphics[width=\textwidth]{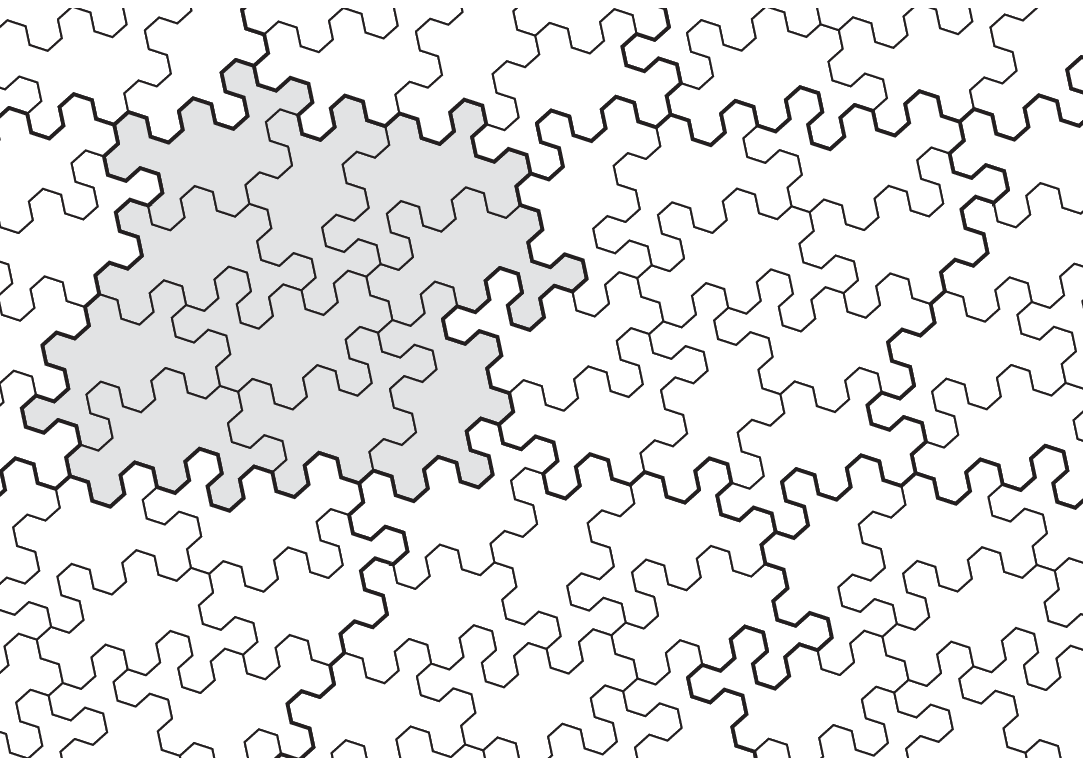}}
\caption{The simplest tessellation by this tile, with isohedral number 10,  the current world-record. The shaded region, found through the brute force enumerations of all possible configurations~\cite{myers}, is a \ddefine{fundamental domain} of the smallest possible size for this tile. }
\end{figure}

Myers has found an amazing variety of strange examples~\cite{myers}--- the two at right in  Figure~\ref{mannmyers} are just the most notable of a great many exotic tiles. Remarkably, the tile at far right has isohedral number 9; the tile just to the left holds the current world-record, with isohedral number 10. But still today, this remains an open question: {\em For all $n$, does there exists a tile with isohedral number at least $n$?}

We might settle the Monotiling Problem through the brute force enumeration of configurations: we can always find out if a tile  does not admit a tessellation, and we can always find out if a tile does admit a periodic tessellation. But it is unknown whether or not there is a procedure to settle this in general, whether or not it is possible to discover whether or not an arbitrary tile does in fact admit a tessellation. Could it be that this problem is {\em undecidable}?

\begin{figure}
\centerline{\includegraphics{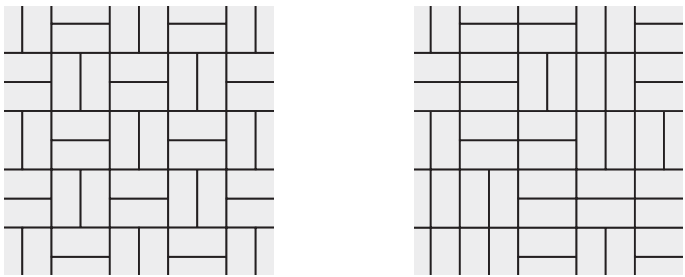}}
\caption{The $2\times 1$ rectangle does admit non-periodic tessellations  ---  for example, one may tessellate the plane with $2\times 2$ squares split horizontally or vertically at random  ---  but admits periodic ones too, and so is {\em not} aperiodic.}\label{dimer}
\end{figure}

This could only be so if there exist tiles that {\em do} admit tessellations of the plane, but {\em do not} admit any periodic tessellations, that is, if there exist \ddefine{aperiodic} tiles. This is far more subtle than merely admitting non-periodic tessellations: as noted in Figure~\ref{dimer}, even the humble $2\times 1$ rectangle manages that simple task! For a tile to be  aperiodic, it somehow would {\em force}  non-periodicity,  somehow  wrecking symmetry at all scales! This is such an astounding property that one might plausibly suppose no such tile could exist, and that the Tiling Problem is in fact decidable. 

The examples we've seen should chasten us:  at the very least, the simple-sounding problem ``How does a tile behave" is astoundingly subtle, and despite much attention, this question has not been conclusively settled! Today, no one can say whether or not the Monotiling Problem is decidable.


\section{Aperiodicity and the Tiling Problem}\label{tiling problem}

More generally, we can ask the 

\vspace{\baselineskip}\noindent{\bf Tiling Problem} or {\bf Domino Problem}: {\em Does a given {\em set} of tiles admit a tessellation of the plane} ---  does there exist  a tessellation formed by copies of tiles in the set?

\vspace{\baselineskip}

This problem arose as an aside in 1961,  as the logician Hao Wang worked on one of the remaining open cases of 
Hilbert's {\em Entscheidungsproblem} (``Is a given first order logical formula satisfiable?") \cite{wang1961}. Wang 
 noted that  the closely related ``Completion Problem" is undecidable and cannot be decided in general by any mechanical procedure whatsoever! 

\label{completion problem} 

\vspace{\baselineskip}\noindent{\bf The Completion Problem:} {\em Does a given a finite collection of tiles admit a tessellation of the plane containing a given ``seed" configuration?}

\vspace{\baselineskip} 

\noindent That is, can we ``complete" the seed configuration to form a tessellation of the entire plane using copies of some or all of the given tiles?

\psfrag{R}[ll][cc][.8]{$t=0$}
\psfrag{S}[ll][cc][.8]{$t=1$}
\psfrag{T}[ll][cc][.8]{$t=2$}
\psfrag{U}[ll][cc][.8]{$t=3$}
\psfrag{V}[ll][cc][.8]{$t=4$}
\psfrag{W}[ll][cc][.8]{$t=5$}
\psfrag{X}[ll][cc][.8]{$t=6$}  
\psfrag{1}[cc][cc][.8]{${\tt 1}$}    
\psfrag{0}[cc][cc][.8]{${\tt 0}$}
\psfrag{A}[cc][cc][.8]{${\tt A}$}
\psfrag{B}[cc][cc][.8]{${\tt B}$}
\psfrag{C}[cc][cc][.8]{${\tt C}$}
\psfrag{A0}[cc][cc][.8]{${\tt A  0}$}
\psfrag{B0}[cc][cc][.8]{${\tt B 0}$}
\psfrag{C0}[cc][cc][.8]{${\tt C 0}$}
\psfrag{A1}[cc][cc][.8]{${\tt A 1}$}
\psfrag{B1}[cc][cc][.8]{${\tt B 1}$}
\psfrag{C1}[cc][cc][.8]{${\tt C 1}$}
\psfrag{pA0}[cc][cc][.6]{$\phi({\tt A  0})$}
\psfrag{pB0}[cc][cc][.6]{$\phi({\tt B 0})$}
\psfrag{pC0}[cc][cc][.6]{$\phi({\tt C 0})$}
\psfrag{pA1}[cc][cc][.6]{$\phi({\tt A 1})$}
\psfrag{pB1}[cc][cc][.6]{$\phi({\tt B 1})$}
\psfrag{pC1}[cc][cc][.6]{$\phi({\tt C 1})$}
\psfrag{A0}[cc][cc][.8]{${\tt A  0}$}
\psfrag{B0}[cc][cc][.8]{${\tt B 0}$}
\psfrag{C0}[cc][cc][.8]{${\tt C 0}$}
\psfrag{A1}[cc][cc][.8]{${\tt A 1}$}
\psfrag{B1}[cc][cc][.8]{${\tt B 1}$}
\psfrag{C1}[cc][cc][.8]{${\tt C 1}$}
\psfrag{H}[cc][cc][.8]{${\tt H}$}
\psfrag{!}[cc][cc][.6]{$\phi({\tt A  0})$}
\psfrag{@}[cc][cc][.6]{$\phi({\tt B 0})$}
\psfrag{#}[cc][cc][.6]{$\phi({\tt C 0})$}
\psfrag{^}[cc][cc][.6]{$\phi({\tt A 1})$}
\psfrag{&}[cc][cc][.6]{$\phi({\tt B 1})$}
\psfrag{*}[cc][cc][.6]{$\phi({\tt C 1})$}
\psfrag{,}[cc][cc][.6]{$\phi({\tt A  0})$}
\psfrag{.}[cc][cc][.6]{$\phi({\tt B 0})$}
\psfrag{/}[cc][cc][.6]{$\phi({\tt C 0})$}
\psfrag{<}[cc][cc][.6]{$\phi({\tt A 1})$}
\psfrag{>}[cc][cc][.6]{$\phi({\tt B 1})$}
\psfrag{?}[cc][cc][.6]{$\phi({\tt C 1})$}

The proof is simple, and rests on encoding Alan Turing's  ``Halting Problem" for machines inside of Wang's Completion Problem for tiles\mp{cross reference to Turing Machine}: For any given Turing machine, Wang produces a set of tiles and seed configuration, in such a way that the seed configuration could be extended to a tessellation by copies of the tiles if and only if the machine never halts. 

Of course the Halting Problem is a touchstone of undecidability: in 1935, through a simple and elegant construction, Turing proved there can be no procedure to decide whether a given Turing machine will halt or not [[cited elsewhere]]; consequently, there can be no procedure to tell whether one of Wang's corresponding sets of tiles, with its seed configuration, can complete a tessellation of the plane.

Wang's construction is easy enough to illustrate by an example: Consider the Turing machine specified by $$
\begin{array}{c|ccc}
 \phi & \tt A & \tt B & \tt C \\ \hline
\tt 0&\tt 0RB &\tt 1LA &\tt 1RB\\
\tt 1&\tt 1RB&\tt 0RC &\tt 0LH
\end{array}
$$

This machine will work on an infinite tape; at each step, each cell is marked $\tt 0$ or $\tt 1$ and the machine will be in a particular state  $\tt A$, $\tt B$ or $\tt C$, reading one particular cell. The transition function $\phi$ determines the action of the machine, depending on its state and the marking it is reading; for example, if the machine is in state $\tt A$ reading $\tt 0$, as in the upper left of the table, the machine will leave a $\tt 0$ in that spot on the tape, move right one cell, and go into state $\tt B$. If it is in state $\tt B$ reading a $\tt 0$, it leaves $\tt 1$ on the tape, moves one cell to the left, and goes into state $\tt A$. There is one special ``halt" state $\tt H$  ---  if the machine enters this state, it can do no more, and the process halts.

\begin{figure}\centerline{{\includegraphics[]{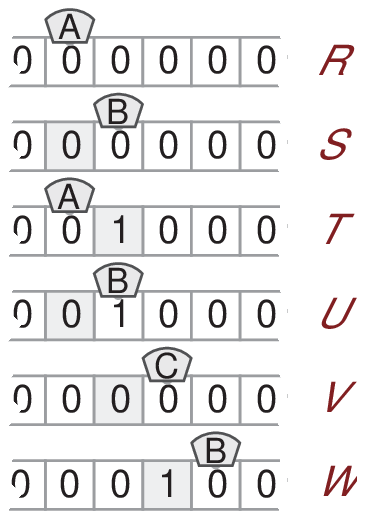}} \hspace{1in}{\includegraphics[]{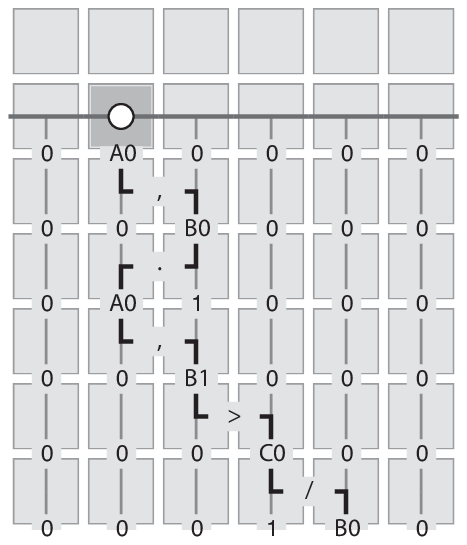}}}\caption{Wang emulates the run of a Turing machine as a tiling problem.}\label{wangrun}\end{figure}

Beginning in state $\tt A$, on a tape marked with all $\tt 0$'s, we can illustrate the first few steps of the run of the machine, at left in Figure~\ref{wangrun}.

The essential point is that this illustration itself satisfies completely local rules: it is composed of pieces that must fit together in a certain manner. We can encode this as a tessellation, shown at right above.

\begin{figure}\centerline{{\includegraphics[scale=.6]{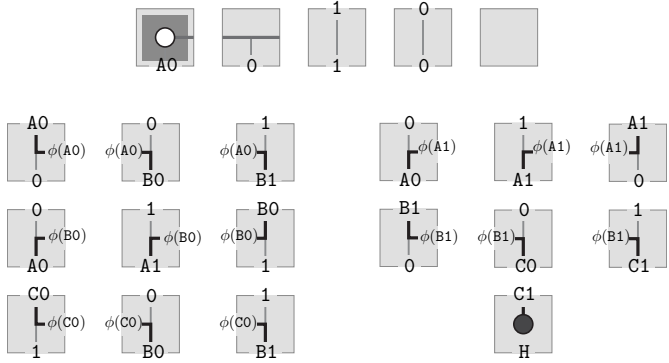}}}\caption{The tiles Wang uses in his Completion Problem, to emulate our sample  machine.}\label{wangMarkings}\end{figure}

With only a little care, we then have a set of tiles, shown in Figure~\ref{wangMarkings}, that {\em can} emulate the machine. It is possible to cover the plane with copies of these tiles, so that labels on adjacent edges match,\footnote{ It is easy enough, if one prefers, to use unmarked tiles that simply are required to fit together: we may convert the labels into geometric jigsaw-like bumps and notches: \includegraphics{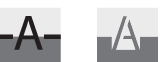}} as shown above.

{\em Must} they emulate this machine? In any tessellation containing the initial ``seed tile", at upper left in Figure~\ref{wangMarkings}, there must be, inductively row by row, a faithful representation of the run of the machine;  this can be completed into a tessellation of the entire plane if and only if the machine never halts  ---  note that the tile at bottom right in Figure~\ref{wangMarkings} corresponds to the machine entering the halt state, and no tile can fit beneath it.
As the Halting Problem is undecidable, so too is the Completion Problem.

(On the other hand, note that it's easy enough to tile in other ways, if we {\em don't} place the seed tile. For example, we could just cover the plane with copies of the filler tile at upper right in Figure~\ref{wangMarkings}.)

This example highlights the deep connection between {\em undecidability} and {\em computational universality}: The celebrated Church-Turing thesis in effect asserts that anything we might mean by computation can be realized by a Turing machine and thus by anything that can emulate a Turing machine. Wang's Completion Problem is undecidable precisely because it has this property, precisely because completing a tessellation from a seed tile is ``computationally universal" and can emulate {\em any} computation (albeit wildly inefficiently!).

As a further aside, Wang posed the Tiling Problem (which is trivial  for the sets constructed above since we may cover the plane with just copies of the blank filler tile). Just as we saw with tessellations by copies of a single tile, Wang noted that if the Domino Problem were in fact undecidable, there must exist sets of tiles that {\em do} admit tessellations, but  {\em none of which} are periodic  ---    if every set of tiles either does not admit a tessellation or admits a periodic tessellation, then we have an algorithm for answering the Domino Problem: enumerate configurations, covering larger and larger disks, until we run out of possibilities ({\tt No}, the tiles do not admit a tessellation), or until we discover a fundamental domain ({\tt Yes}, the tiles admit a tessellation). 

Wang reasonably conjectured no such \ddefine{aperiodic set of tiles} could exist  ---    after all, somehow, just by local rules, symmetry would have to broken at all scales  ---    but within a few years Robert Berger, and then Raphael Robinson gave subtle proofs that the Domino Problem is undecidable,  along the way producing aperiodic sets of tiles~\cite{berger, rob}.

 Berger's proof that his set is aperiodic rests on their forcing the emergence of a particular kind of non-periodic, hierarchical structure; a few years later Robinson gave a greatly simplified construction which we momentarily describe. Both authors used these forced hierarchical structure as a kind of scaffolding in their demonstrations of the undecidability of the Tiling Problem, as we will soon sketch in Section~\ref{domino}. The construction has been widely imitated and generalized, producing many other aperiodic sets of tiles, as we discuss in Section~\ref{aht}.
 
\subsection{Robinson's aperiodic set of tiles}

\begin{figure}\centerline{\includegraphics{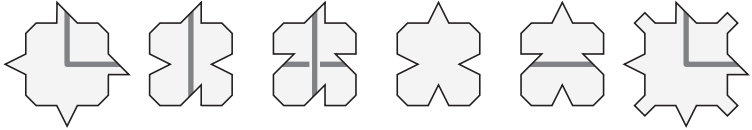}}\caption{The Robinson tiles are aperiodic.}\label{fig_robinson_1}\end{figure}
  
  Robinson's six tiles are shown in Figure~\ref{fig_robinson_1}. Let us show that: 
  
\vspace{\baselineskip}\noindent{\bf Theorem }{ \em The Robinson tiles are an aperiodic set.}

\vspace{\baselineskip}

 \begin{proof} We have two tasks: we must show that the Robinson tiles {\em can tile} the plane, but also that they {\em cannot tile periodically}.
We first  notice  that the last tile in Figure~\ref{fig_robinson_1} is ``cornered" and the other five are ``cornerless."  The first and last ``cross" tiles have a ``bump" pointing out of each edge; the other four ``passing" tiles have one ``bump" and three ``nicks". The two cross tiles have an ``elbow", shown as pointing up and to the right.

\begin{figure}\centerline{\includegraphics{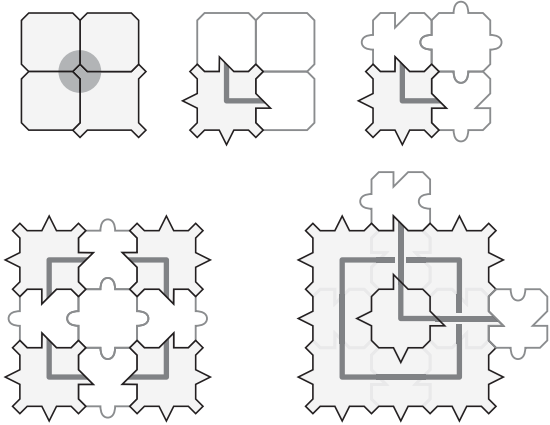}}\caption{Robinson's proof rests on showing that tiles form blocks that themselves fit together as the tiles do, on a larger scale, forming still larger blocks, {\em ad infinitum}.}\label{robinson_0}\end{figure}

 In any configuration of these tiles, at every corner, three cornerless tiles must meet one cornered tile, as shown at top left in Figure~\ref{robinson_0}.


What tiles can be filled in the positions adjacent to the elbow of a cornered cross tile, as indicated in the top middle of Figure~\ref{robinson_0}?
We can only have two passing tiles, as shown at top right, and thus a cornerless cross adjacent to those.  
We quickly see that every cornered cross can only belong to the $3\times 3$ configuration shown at bottom left in Figure~\ref{robinson_0} ---  though we do have a choice for the orientation of the central elbow.


In turn, at the ends of this larger central elbow, only two passing tiles can fit, as shown at bottom right in Figure~\ref{robinson_0}, and hence a chain of passing tiles, forcing the placement of a cornerless cross, at the upper right of the figure above. In just the same way that the original tiles could only form $3\times 3$ blocks, these $3\times 3$ blocks can only fit together to form a $7\times 7$ block, 
 these $7\times 7$ blocks can only fit together $15\times 15$ blocks, and these into $31\times 31$ blocks, and so on. In Figure~\ref{fig_robinson_7} we can see a nested sequence of such blocks.

\begin{figure}\centerline{\includegraphics{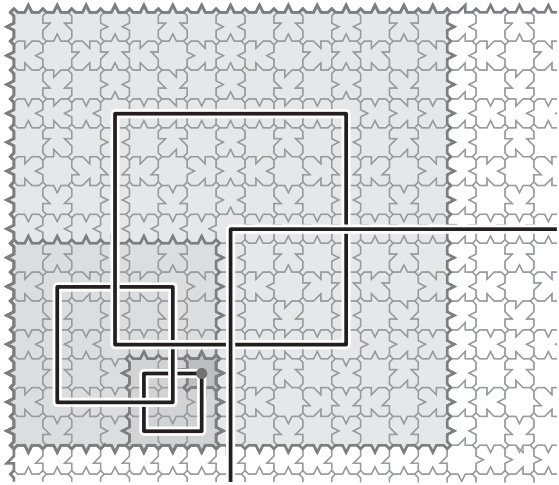}}\caption{Constructing one of uncountably many tessellations admitted by the Robinson tiles, as a countable sequence of ever more momentous choices.}\label{fig_robinson_7}\end{figure}

{\em Do the Robinson tiles admit a tessellation?} Let us explicitly construct a tessellation: fix a crossed corner tile. We may place a $3\times 3$ block to the southwest; this in turn may be in the northwest of a $7\times 7$ block, which will be in the northeast of a $15\times 15$ block, which will in its turn be in the southeast of a $31\times 31$ block, and so on. Continuing in this way, we have a well-defined tessellation of the plane: every point in the plane will lie within a specified well-defined tile.

In fact, there are uncountably many distinct tessellations by the Robinson tiles! We had an infinite sequence of choices in the construction: at each stage, the next largest block might have lain in any of four positions relative to what we had already placed. 
There are two subtleties: not every sequence actually will produce a configuration that covers the entire plane (for example, if we always placed the next block to the southwest). But still uncountably many do. 
Different sequences of choices might have given congruent tessellations, but as any given tessellation would have only been congruent to countably many others, the Robinson tiles still admit uncountably many non-congruent tessellations of the plane.

\begin{figure}\centerline{\includegraphics{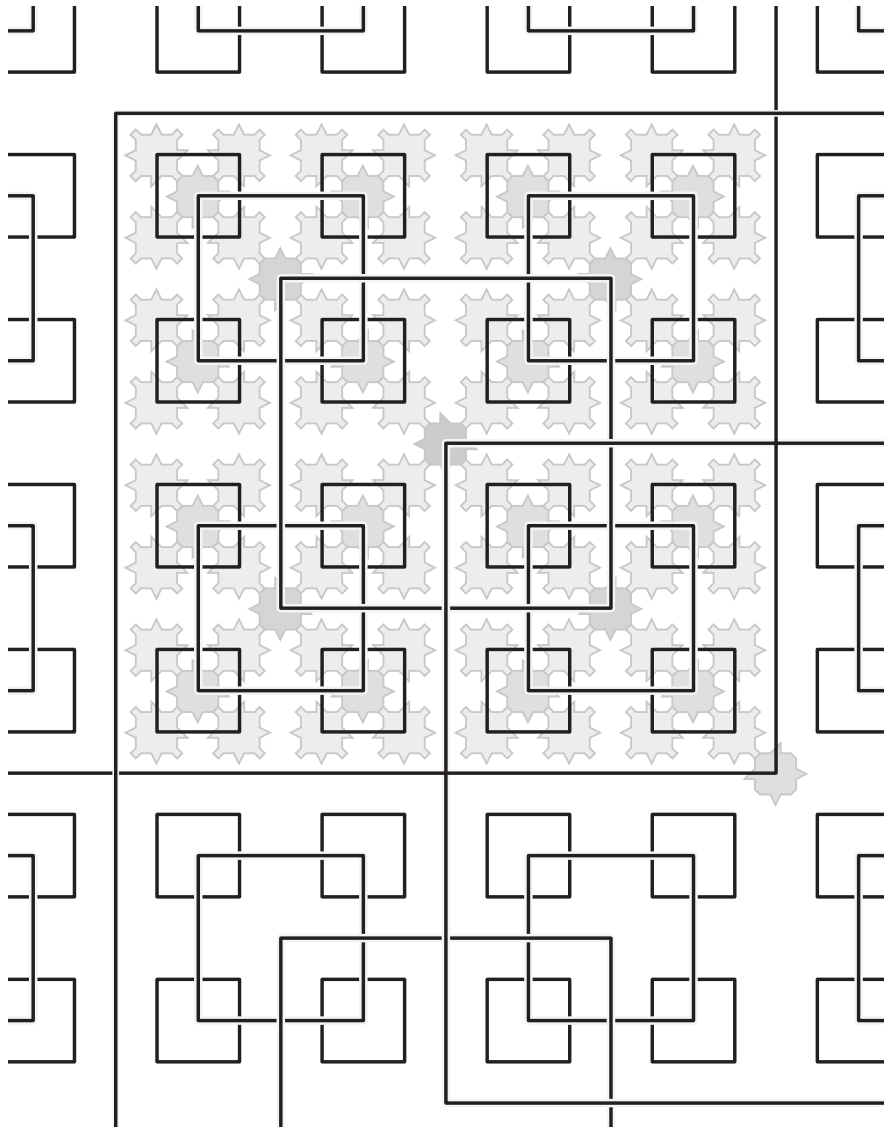}}\caption{A heirarchy of squares: the Robinson tiles do not admit periodic tessellations. Note the complementary systems of overpasses and branches.}\label{fig_robinson_6}\end{figure}

{\em Can the Robinson tiles admit a periodic tessellation?} As shown in Figure~\ref{fig_robinson_6}, a cross tile in the center of an $(2^n-1)\times (2^n-1)$ block lies at the vertex of a marked $2^{n-1}$ by $2^{n-1}$ square; squares of the same size do not overlap. 

Suppose there were a periodic tessellation admitted by the Robinson tiles. 
Then translation by some vector $\vect v$ would leave the tessellation invariant. (That is, if you were to turn your back while I shift by $\vect v$, when you again take a look you would be unable to tell that anything had happened.) 

But this is impossible: there is some size square $S$ so large that $\vect v$ would be unable to take $S$ to another square of the same size, and so shifting by $\vect v$ could not leave the tessellation invariant.
\end{proof}

\subsection{Aperiodic hierarchical tessellations}\label{aht}

Berger's great discovery, that hierarchical structure can be forced through local rules, led to  many notable and beautiful examples. The Robinson tiles we just saw are quite typical in many respects: from some hierarchical structure we {\em wish} to form, in that case, a hierarchy of larger and larger squares, as in Figure~\ref{fig_robinson_6}, we then find a way of marking tiles that {\em forces} that structure to emerge. Frankly, in all the simple examples, there is a touch of magic to this!

Once we somehow have actually found such tiles, the proof that they actually do what we wish, that they actually force such a hierarchical structure,  tends to follow familiar lines: we show they must form larger blocks, which themselves behave just as the original tiles do, and so must themselves form blocks that are larger still, which... {\em ad infinitum}. 

But what hierarchical structures are there in the first place, for us to try to enforce?
Before we describe a few more examples, let us pause for a moment to discuss {\em substitution tessellations}

\subsection{Rep-tiles and substitutions}

\begin{figure}\centerline{\includegraphics[width=.6\textwidth]{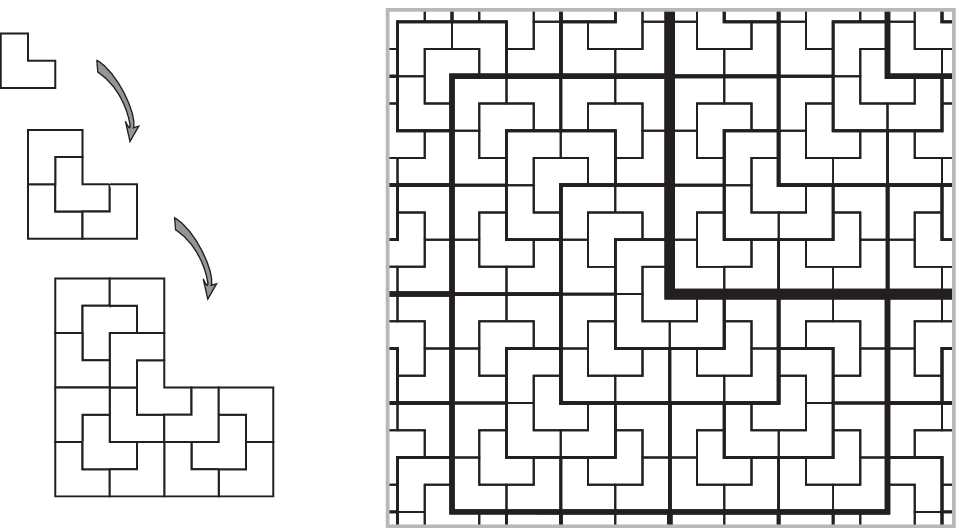}}\caption{A rep-tile, the ``chair", and a corresponding hierarchical tessellation.}\label{L}\end{figure}

An old puzzle asks how an farmer can divide his L-shaped piece of land into four plots for his for jealous sons; each plot must be just the same size and shape as the others. The answer is shown in Figure~\ref{L}: a large L can be subdivided into four small L's. A tile that can be subdivided into smaller copies of itself is often called a ``rep-tile", which is kind of a pun in English (it is made of replicas of itself). Rep-tiles readily can be used to form hierarchical structures; just as in Figure~\ref{fig_robinson_7}, we can assemble our tiles into larger copies, which we then assemble into larger copies still, and so on. With care we can assemble a tessellation of the entire plane, one that automatically comes with a hierarchical structure. 

\begin{figure}\centerline{\includegraphics[width=.8\textwidth]{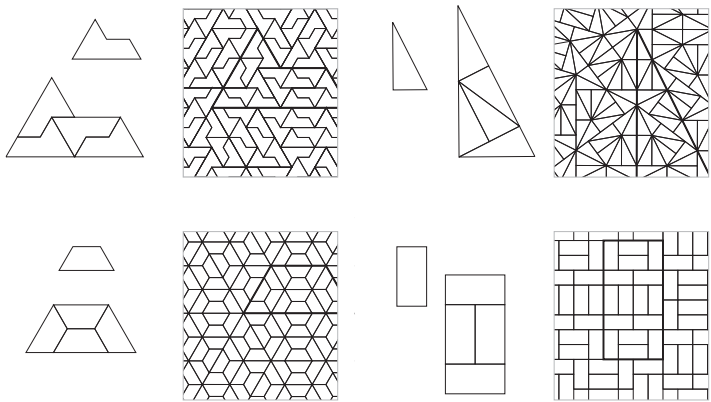}}\caption{A variety of rep-tiles and corresponding tessellations, clockwise from top left: the sphinx,  Conway's pinwheel, the dimer and the half hex. }\label{rep}\end{figure}

\begin{figure}\centerline{\includegraphics[width=.8\textwidth]{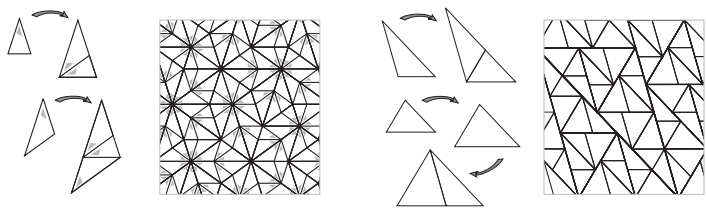}}\caption{Two triangle substitutions, at left ``gold" and at right ``silver", with corresponding tessellations.}\label{tri}\end{figure}

A few more rep-tiles are shown in Figure~\ref{rep}, with portions of a hierarchical tiling. More generally,  we might have two or more different kinds of tile, subdivided into copies of themselves in some more complex manner. For example, the two ``golden" triangles at left in Figure~\ref{tri} can each be subdivided into copies of themselves, as can the ``silver" triangles at right.

\begin{figure}\centerline{\includegraphics[width=.85\textwidth]{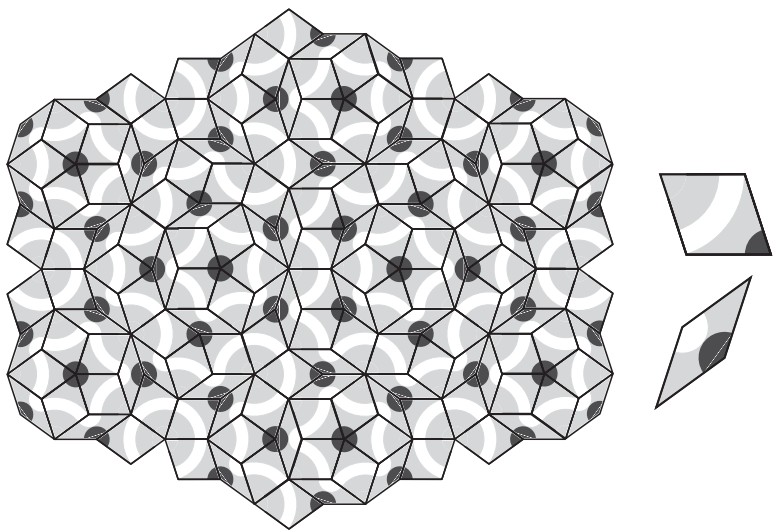}}
\caption{A configuration of Penrose rhombs; note the sinuous ribbons, the ever irregular arrangement of black dots, and the appearance of the rhombs as cubes in perspective ---  all with nested five-fold symmetries.}\label{penrose_rhombs}
\end{figure}
 
 Many examples of these kinds of substitutions on tiles are known ---Ludwig Danzer in particular has found many infinite families with very strange properties~\cite{tilingzoo}--- but the fact is, no one really understands what is possible, or what is not, or why. There is room for many more discoveries!
 
\subsection{Aperiodic tiles}

Now the rep-tiles we've just seen are not themselves {\em aperiodic}. The unmarked tiles can form a great variety of tessellations, some of which are periodic and some of which, like the hierarchical tilings we've illustrated, are non-periodic. The tiles do not {\em force} this structure. Which hierarchical structures {\em can} be forced?

Once Berger pointed the way, a number of notable, elegant examples were found. We've just seen Robinson's, after which Roger Penrose soon followed, forcing hierarchical structure based on subdivisions of the pentagon. His rhombs (Figure~\ref{penrose_rhombs}) and his kites and darts~\cite{gardner} are probably the most famous of all aperiodic sets of tiles. In the late 1970's Robert Amman gave several remarkable examples, finally published in {\em Tilings and Patterns}~\cite{grsh}--- still today, his tiles are among the very simplest known. Ludwig Danzer, Joshua Socolar, Charles Radin, myself and a few others have given additional sets of this general kind. But 
all together, remarkably few sets were known, and all seemed like art, or magic.

In 1988, Sharir Mozes gave the first general construction, for systems based on special subdivisions of the rectangle. For any subdivision of a particular kind, he was able to produce new tiles that faithfully forced the corresponding hierarchical structure to emerge. 

At this point, may I modestly report that the complete solution--- quite indebted to all these examples--- appeared in~\cite{mrst}? Given virtually any substitution ---  and there are certainly infinite familes of examples ---  we can produce a set of tiles to force it. But still, this construction is so general that in practice it produces large, unwieldy sets of tiles. There is still room for novel, simple, magical constructions, and from time to time a new one is found.

In particular, towards the end of Section~\ref{whattiles}, we glossed over a particular famous question: {\em Does there exist an ``einstein", an aperiodic monotile, that is, a single tile that admits only non-periodic tessellations?}

To some extent, we've found that this becomes almost a semantic, or perhaps sociological question--- what do we mean, really? Are we allowed to mark our tiles? Can we allow more complex matching rules than just fitting the tiles together? What if we allow the tiles to overlap? Who decides, and how, what the rules should be?

Generalizing very slightly, Penrose gave a reasonable example~\cite{penrosemono} and just recently Joan Taylor and Socolar gave another one, and David Fletcher yet another. But if we require unmarked tiles, simply fitting together, as shown throughout Section~\ref{whattiles}, then we do not know the answer, whether or not there is an ``einstein".

More critically, what other kinds of aperiodic sets of tiles are possible? All of the examples we have seen force hierarchical structure, but by no means is that the only possibility. Space limitations do not allow discussion here, but there is a literature on ``quasiperiodic" tiles, sets of tiles which force a completely different kind of non-periodic structure\cite{debruijn,le}, and in 1996, Jarkko Kari gave a completely different means of forming aperiodic tiles~\cite{kari}. But that's it. These three methods are just a shadow of what is possible. Precisely because of the undecidability of the Tiling Problem, discussed in a moment, we know that we can  never fully answer this question: {\em What are other means of constructing aperiodic sets of tiles?}  Mysteries abound.

\section{The Undecidability of the Tiling Problem}\label{domino}

And this will always be so. As we mentioned in Section~\ref{tiling problem}, Berger {\em proved} that we can never find a general solution to the Tiling Problem, we can never fully understand how to determine whether or not any given set of tiles admits a tiling. 

Just as Robinson simplified Berger's initial construction of an aperiodic set of tiles, he also streamlined Berger's proof that the Tiling Problem is undecidable, and we can at least get the idea here. Recall that Wang's Completion problem was proven undecidable by emulating the action of a Turing Machine, but required the placement of an initial ``seed" configuration. Wang himself could not see how to {\em force} this placement and still carry out the emulation. Berger's idea was instead to begin by constructing a set of tiles that can only admit a non-periodic, hierarchical structure. Within this structure, we hang increasingly large domains in which we can emulate finite, but longer and longer runs of a given Turing machine.  All of these domains can be tiled, fully, if and only if the machine itself does not halt.

\begin{figure}\centerline{\includegraphics{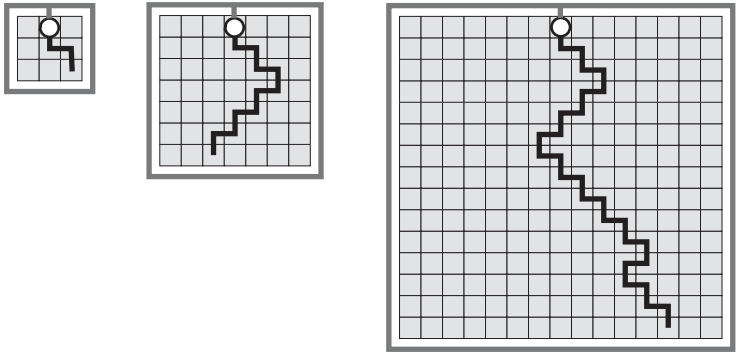}}\caption{Berger's basic idea for the proof the undecidability of the Tiling Problem: emulate Wang's proof, on larger and larger domains. (Compare to Figure~\ref{wangrun}.)}\label{berger1}\end{figure}

In Figure~\ref{berger1}, we sketch the basic idea.  Once the corners of each domain are placed, as our ``seeds", the remaining tiles will be forced, emulating the behavior of our given Turing machine for some limited period of time. We will be able to complete domains of arbitrary size ---  and thus a tessellation of the plane ---  if and only if the Turing machine does not halt. As this is undecidable, so too will be the Tiling Problem. 

\begin{figure}\centerline{\includegraphics[width = .85\textwidth]{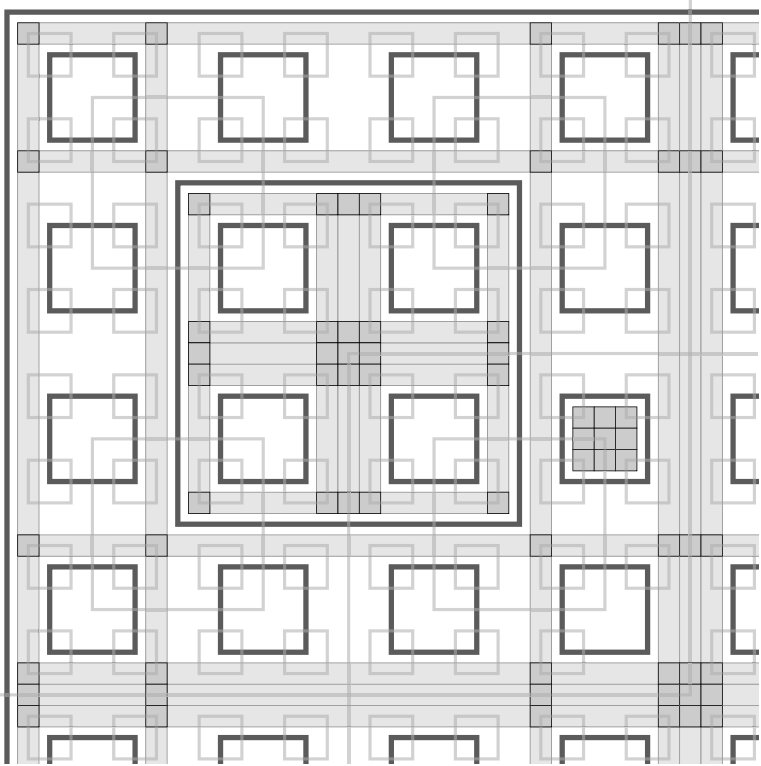}}\caption{Here is the heart of Robinson's proof that the Tiling Problem is undecidable. Onto his aperiodic set of tiles, we paint a variety of new markings that encode domains and the Turing machines that are emulated within them. The domains will be nested as shown: $3\times 3$ domains within the $4\times 4$ blocks of the second level of the original hierarchy, then $5\times 5$ domains within the $16\times 16$ blocks of the fourth level, $9\times 9$ domains within $64\times 64$ blocks, etc. We can complete a tessellation by the tiles we thus produce if and only if we can emulate arbitrarily long runs of the machine --- that is, if and only if the machine does not halt.}\label{robinsonProof}\end{figure}

But now, how can we ensure these arbitrarily large domains are themselves present? In Figure~\ref{robinsonProof} we illustrate Robinson's solution: let us ignore every other level of squares in the hierarchy forced by the Robinson tiles; a sequence of nested squares remain ---  these form our domains. 

\psfrag{a}[cc][cc]{$a$}

\begin{figure}\centerline{\includegraphics[]{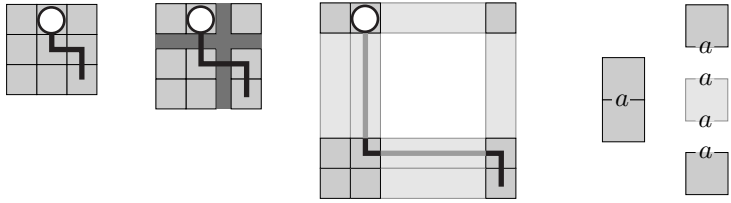}}\caption{Our domains will be nested, and as at left larger domains will have to be attenuated in order to ``flow" around smaller ones. We allow this by introducing new tiles. For example, if  an edge  marked $a$ in our original simulation were to be attenuated to allow room for a smaller domain, we introduce a tile marked as shown at right.
}\label{robinson_10}\end{figure}

We need one final trick: smaller domains lie within larger ones, getting in the way of our emulation. All we must do is to flow around them, attenuating our simulation, as sketched in Figure~\ref{robinson_10}. Where two tiles would once have met, there is now a string of what we might call ``edge tiles"  ---  tiles which emulate the markings on our original edges.

Returning to Figure~\ref{robinsonProof}, we see the first three levels of nested domains, attenuated to flow around those that are smaller, leaving $3\times 3$, then $5\times 5$, then $9\times 9$ cells for the runs of a machine, using squares of the 2nd, 4th and 6th levels of Robinson's hierarchical structure. In general, from the $(2n)$th level, we obtain a domain of $2^n+1$ cells width and height.

With the idea in hand, the rest is fairly technical but relatively straightforward. We overlay the Robinson tiles with additional markings to encode this additional structure, as well as with markings similar to those in 
 Figure~\ref{wangMarkings}: We now can reliably emulate finite runs of the machine. For a given machine, the corresponding set of tiles admits a tiling if and only if the machine runs for arbitrarily long time  ---  that is, if and only if the machine does not ever halt.

We know so little, and there is so much to explore. This topic is rich with questions that need little formal training --- but much imagination!--- to explore, and I hope many readers take up the challenge!

I'd like to conclude with a final image, of a single shape of tile that can form remarkable tessellations with scaled copies of itself.  These tiles, when drawn to infinite detail, are all equivalent up to scaling, and admit myriad tessellations that cover {\em almost} the entire plane ---  the structure of these missing points in such tilings is a rich and fascinating open area to explore!

\begin{figure}\centerline{\includegraphics[width=\textwidth]{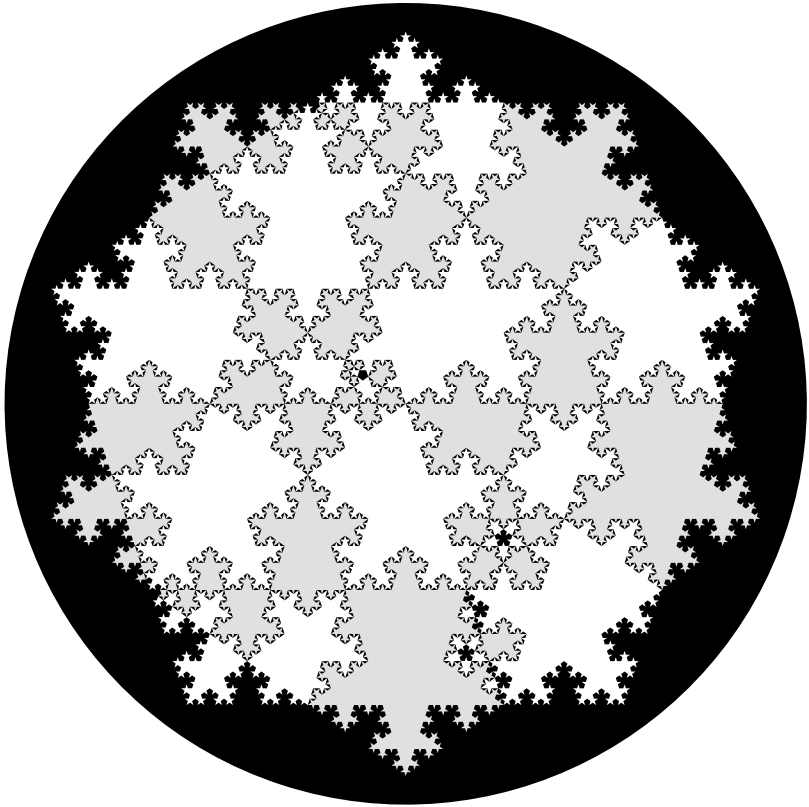}}\caption{What can copies  this tile form? What structure is there in the holes in tessellations by this tile? What other self-similar tiles can you find?}
\label{pentakoch}\end{figure}

\appendix\section{The Euclidean Tiling Theorem}

In Section~\ref{whattiles} we promised a statement of the powerful Euclidean Tiling Theorem; it's well worth proving, too, since the techniques are so elementary, yet broad and useful. The key lemma applies, in fact, to tilings of any surface!

Before we begin, we note there's no limit on how many {\em sides} a polygon can have and still admit a tiling of the plane, if we allow adjacent pairs of tiles to meet along as many successive edges as we please:

\vspace{\baselineskip}
\centerline{\includegraphics[]{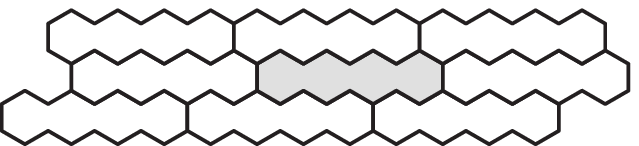}}

But there is a strict constraint on the number of  neighbors a tile might have, or the number of tiles that can meet at a vertex, at least on average, at least if we require all our tiles to be about the same size. For example, it's perfectly possible to tessellate with tiles meeting three-to-a-vertex, each having seven neighbors, but these heptagons will have to be increasingly distorted in order to fit together. On the other hand, if all of our heptagons are about the same shape and size, or even are exactly the same shape and size as in all the examples we've seen so far, then this is impossible. 

In essence, we will define a way to measure a kind of combinatorial imbalance--- do our tiles have too many neighbors, or too few? The theorem will state that this imbalance must average out to nearly zero over very large configurations.

We make this precise: Consider a tessellation by tiles that are ``reasonably round": the tiles are topological disks and if a configuration roughly covers a large disk, then the number of tiles within it is roughly proportional to the area of the disk and the number of edges on its boundary is roughly proportional to the circumference of the disk. More precisely, a tessellation is ``reasonably round" if there are constants $c_1$ and $c_2$ so that in any disk of radius $R$, there are at least $c_1 R^2$ tiles and the boundary of the configuration of these tiles has no more than $c_2 R$ edges. 

This can be guaranteed in many ways, and a complete discussion can be found in Chapter 3 of Gr\"unbaum and Shephard's masterpiece {\em Tilings and Patterns}~\cite{grsh}. For example, we can simply require that all of our tiles can inscribe, and in turn be inscribed in a pair of particularly sized disks, and that they all have a bounded number of neighbors. But we can loosen even these restrictions somewhat and still maintain the key property of being ``reasonably round." 

In a Euclidean tessellation, there can be tiles with a great many neighbors, or a tremendous number of tiles meeting at a given vertex, but these have to be balanced out by tiles with just  a few neighbors, or a small number of tiles meeting at some other vertex. This balancing condition is very precise, and fairly easy to state. In order to do so, let us contrive a way to measure how far from ``balanced" a given tile is.

Suppose a given tile meets $n$ others in the tiling (if, somehow, a tile meets another more than once, we just count with multiplicity), and that reading around, at the first vertex of the tile, $q_1$ tiles meet, then $q_2$, etc., until at the last vertex, $q_n$ tiles meet. Let's call $2\pi/q_1, 2\pi/q_2, \ldots$ ``valence angles" of the polygon; then the ``imbalance" of the tile is $\Sigma -  (n-2)\pi$, where $\Sigma$ is the sum of all the valence angles. This looks strange, but this has a simple interpretation. Remember that the sum of the {\em vertex} angles of any $n$-sided polygon in the Euclidean plane is exactly $(n-2)\pi$. This imbalance measures how far we are from being able to straighten out all the valence angles into vertex angles.\footnote{We may regard this imbalance as a kind of combinatorial curvature which exactly measures the total curvature the tiles would have to have if were to modify them so the valence angles were their actual vertex angles.  For example, we might modify the tessellation by octagons and squares, so each vertex angle really is $120^{\circ}$.  An octagon with $120^{\circ}$ vertex angles must have negative curvature, totalling $-2\pi/3$ over its surface. Similarly, a square with $120^{\circ}$ vertex angles must have positive curvature ---  one way to construct such a square is to take one face of the spherical cube, with total curvature one sixth the total curvature of the sphere, or $4\pi
\div 6 = 2\pi/3$. In essence all that follows here can be summed up in these terms: The total curvature of  tessellation by tiles with curvature must, in the large, approach zero if we are to be able to flatten it out into a Euclidean tessellation without too much distortion.

But how can we use $K$ as a stand in for curvature? 
The celebrated Gauss-Bonet theorem relates the total curvature of a disk-like region to the net amount of turning we experience as we traverse its boundary.  If the region has zero-curvature, this net turning will be exactly $2\pi$  ---  try it on any flat disk-like region in the plane! But around a positively curved region, this turning will be less than $2\pi$ and if the total curvature is negative, the turning will be greater than $2\pi$. The term ``holonomy" is shorthand for this difference: it is just $2\pi$ minus the amount we turn as we traverse the boundary of a region. It is remarkable that the total curvature, $\int_S \kappa$ (that is, the total of all the Gaussian curvature $\kappa$ across $S$), of a disk-like region $S$ is exactly the same as the holonomy around its boundary. 

Now what happens as we walk around a polygonal disk with $n$ straight sides, but possibly having some curvature across its interior?
If the interior angle at a given vertex is $\alpha$, then we turn $\pi - \alpha$ at that vertex. If $\Sigma$ is the sum of all these interior angles, the total turning is $\pi n - \Sigma$ and our holonomy, and hence total curvature, is $2\pi - (\pi n - \Sigma) = \Sigma - (n-2)\pi$ ---  exactly the same as $K$ if our ``valence angles" really were our actual interior angles.}

 \vspace{\baselineskip}
\centerline{\includegraphics{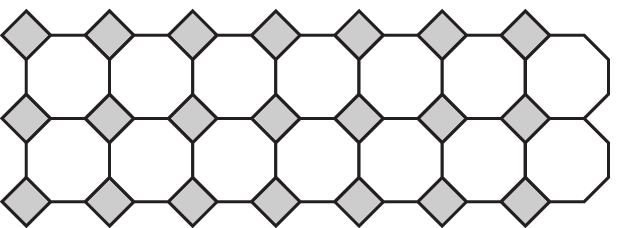}}
 \vspace{\baselineskip}

For example, in the tessellation by octagons and squares, each vertex has valence three, and all the valence angles are $2\pi/3$.  Each octagon has an imbalance of $8\times 2\pi/3 - (8-2)\pi  = -2/3\pi$ and each square has an imbalance of $4\times 2\pi/3 - (4-2)\pi = +2/3\pi$. In the configuration pictured here the imbalances don't quite cancel out ---  there are more squares than octagons ---  but in a very large region, there will be about the same number of octagons as squares, at least relative to the total number of tiles, and the net imbalance, relative to the area of the configuration, will approach zero. This next theorem provides strong topological restrictions on Euclidean tessellations:

\vspace{\baselineskip} 
\noindent
{\bf The Euclidean Tessellation Theorem: } {\em For any Euclidean tessellation by reasonably round tiles, considering configurations $C(r)$ of tiles meeting larger and larger disks of radius $r$, the average imbalance of the tiles in $C(r)$ must approach zero.}

\vspace{\baselineskip} 

The following lemma contains the essential truth of the matter: let us define the imbalance of a {\em configuration} in the same way we did for individual tiles. A configuration is itself just a large polygon  subdivided into tiles. Ignoring the tiles themselves, suppose the configuration has $n$ vertices, with valences $q_1, \ldots q_n$; then define the imbalance to be $\Sigma - (n-2)\pi$, where $\Sigma$ is the sum of the valence angles $2\pi/q_1,  \ldots, 2\pi/q_n$. We'll pause for this deceptively powerful lemma  ---  which holds even for tessellations in other spaces
such as the sphere or hyperbolic plane  ---  and then use it in the proof of our theorem: 

\vspace{\baselineskip} 
\noindent
{\bf The Tiling Lemma: } {\em The imbalance of any configuration is exactly equal to the sum of the imbalances of the tiles within it. Consequently, if a configuration has imbalance $K$ and $n$ vertices on its boundary; then $|K-2\pi|<\pi n$}

\vspace{\baselineskip} 

One route is to prove it inductively, chopping our configuration into smaller ones for which we presume the lemma holds. We'll prove the lemma, though, using Euler's celebrated theorem for planar maps, giving  us a stringent restriction on the numbers of vertices ($v$), edges ($e$) and tiles ($f$) in any bounded configuration:
 $$v-e+f=1$$ This does require us to continue to make two assumptions we've already been using: that 
  our configuration and tiles are  topological disks. 
  
  Try it! For example in the configuration of octagons and squares we saw earlier, there are ninety-one vertices, one hundred twenty-five edges and thirty-five tiles,  and indeed, $91 - 125 + 35 = 1$! Euler's Theorem generalizes readily and plays  a large role in shaping a wide range of mathematical phenomena, such as vector fields or the topology of surfaces. Here we can get control of the imbalance of configurations, and thus exactly how many neighbors, on the whole, can a tile have in any tessellation.
  
   All we really have to do is make a careful count of all the valence angles in all the tiles, applying Euler's Theorem at just the right moment. With that in mind, let's fix a configuration $C$, and  let $v,e$ and $f$ be the number of vertices, edges and tiles in $C$. Let $\Sigma$ be the sum of the valence angles of $C$ and $n$ be the number of vertices on the boundary of $C$.  

For a tile $t\in C$, let $K_t$ be the tile's imbalance, $n_t$ be the the number of edges it has, and $\Sigma_t$ be the sum of its valence angles. Then summing this imbalance over all the tiles in $C$, we have
  $$\sum_{t\in C} K_t  = \sum_{t\in C} \Sigma_t -  (n_t-2)\pi $$

Now $\sum_{t\in C} \Sigma_t$ is the sum of all the valence angles in the entire tiling. Each of the $v-n$ vertices in the interior of $C$ contributes exactly $2\pi$ to this sum, and all the rest contribute a total of $\Sigma$. Therefore $\sum_{t\in C} \Sigma_t =  2\pi(v-n) + \Sigma$.

Next, $\sum_{t\in C} n_t$ counts the number of edges of each tile ---  this counts each of the $e-n$ edges in the interior of $C$ twice (since each such edge meets two tiles) and each of the $n$ edges on the boundary of $C$ just once. So $\sum_{t\in C} n_t = 2(e-n)+n = 2e - n$ and so $\sum_{t\in C} n_t \pi = 2\pi e - \pi n$. 

Lastly $\sum_{t\in C} 2\pi= 2\pi f$, giving us a grand total of $$\sum_{t\in C} K_t   = 2\pi(v-n) + \Sigma - (2\pi e - \pi n) + 2\pi f = \Sigma - \pi n + 2\pi (v+f-e) = \Sigma - (n-2)\pi$$ The sum of the imbalances of the tiles is indeed the imbalance of the configuration! 

Now we restate this, obtaining our bound on the imbalance $K$: 
Each valence angle is more than $0$ but less than $2\pi$ and consequently the sum $\Sigma$ of all the valence angles is more than $0$ but less than $2\pi n$. Remembering that $K=\Sigma - (n-2)\pi$ we have $0-(n-2)\pi < K  < 2\pi n - (n-2)\pi$ or $$-\pi n < K - 2\pi < \pi n$$ 
In other words, the imbalance of a configuration is very strictly bounded by the length of its boundary. \qed

It is very nice that this observation holds no matter what space we are tiling: it is just as true on the sphere or in the hyperbolic plane! What is different in the Euclidean plane is that the area of a disk grows much faster than the radius of a disk. Our Euclidean Tessellation Theorem follows immediately from the lemma: 

If we take configurations $C(r)$ of reasonably round tiles meeting disks of larger and larger radius $r$,  by the definition of reasonably round there will be constants $c_1$ and $c_2$ so that  the number of tiles $N_r$ in $C(r)$ will be at least $c_1 r^2$ and the number of vertices on the boundary of $C(r)$ will be less than some $c_2 r$. By our lemma, the imbalance $K_r$ of $C(r)$ will satisfy $-(c_2 r\pi - 2\pi) <K< c_2 r \pi + 2\pi$ and $|K|< c_2 r \pi - 2\pi$;  thus $|K_r|/N_r<  (c_2 r \pi - 2\pi)/(c_1 r^2)$, which tends to zero as $r$ increases.\qed

To sum this all up,  we have shown that in a Euclidean tessellation by ``reasonably round" tiles, on average, over larger and larger regions, the average imbalance must be closer and closer to zero. This has many immediate corollaries. For example:

\vspace{\baselineskip}
\noindent
 {\bf Corollary: } {\em No tile admits a Euclidean tessellation in such a way  that every copy of the tile has seven or more neighbors.}
 
 \vspace{\baselineskip}
 
 Any tessellation by copies of a single tile is reasonably round automatically ---  all the tiles are exactly congruent. If every tile has seven neighbors, then every tile is negatively imbalanced, and so the average cannot approach zero. The number of sides does matter if a tile is {\em convex}: in a tessellation by {\em convex} tiles, every tile will have at least as many neighbors as it has sides. We thus know, for example:
 
 \vspace{\baselineskip} \noindent
 {\bf Corollary: } {\em No convex $n$-gon, $n\geq 7$,  admits a Euclidean tessellation.}



\begin{thebibliography}{6}

\bibitem{berger} {R. Berger}, {\em The undecidability of the Domino Problem},
Memoirs Am. Math. Soc.  {\bf 66} (1966).

\bibitem{sot} J.H. Conway, H. Burgiel and C. Goodman-Strauss {\em The Symmetries of Things}, A.K. Peters (2008).

\bibitem{debruijn} N.G. de Bruijn, {\em Algebraic theory of Penrose's
non-periodic tilings}, Nederl. Akad. Wentensch Proc. Ser. A   {\bf 84} (1981),
 39-66.

\bibitem{tilingzoo} D. Frettl\"oh and E. Harriss {Tilings Encyclopedia} {\em http://tilings.math.uni-bielefeld,de/tilings/index}

\bibitem{gardner} M. Gardner, 
{\em Extraordinary nonperiodic tiling that enriches the theory of tilings},
 Scientific American  {\bf 236} (1977), 110-121.

\bibitem{mrst} {C. Goodman-Strauss}, {\em Matching rules and substitution
tilings},  Annals of Math. {\bf 147} (1998), 181-223.



\bibitem{grsh}{B. Gr\"unbaum} and {G.C. Shepherd}, {\em Tilings and patterns},
W.H. Freeman and Co.   (1987).

\bibitem{kari} {J. Kari}, {\em A small aperiodic set of Wang tiles},
Discrete Mathematics {\bf 160} (1996), 259-264.

\bibitem{klarner} {\em Mathematical Recreations: A Collection in Honor of Martin Gardner}, D. Klarner (ed.), Dover Publications (1998).

\bibitem{le} {T.T.Q. Le}, {\em Local rules for quasiperiodic 
tilings}, Proc. of NATO ASI Series C 489, ``The Mathematics of
Long Range Aperiodic order'', (ed. R.V. Moody), Kluwer (1997), pp. 331-366.

\bibitem{mannHeesch} C. Mann {\em Heesch's tiling problem}, Amer Math Monthly {\bf 111} (2004), 509-517.

\bibitem{mann07}  Personal communication.

\bibitem{mozes} {S. Mozes}, {\em Tilings, substitution systems and dynamical
systems generated by them}, J. D'Analyse Math.  {\bf 53} (1989), 139-186.


\bibitem{myers} J. Myers {\em Polyomino, polyhex and polyiamond tiling}, http://www.srcf.ucam.org/$\sim$jsm28/tiling/

\bibitem{pegg} E. Pegg, {\em Pentagon Tilings}, a Wolfram Demonstrations Project, {\em http://demonstrations.wolfram.com/PentagonTilings/}

\bibitem{penrosemono} R. Penrose {\em Remarks on Tiling: details of a
$(1+\epsilon+\epsilon^2)$-aperiodic set}, The mathematics long range
aperiodic order, NATO Adv. Sci. Inst. Ser. C. Math. Phys. Sci. 489 (1997),
467-497.


\bibitem{rob} {R. Robinson}, {\em Undecidability and nonperiodicity of tilings
in the plane}, Inv. Math.  {\bf 12} (1971), 177-209.

\bibitem{schattschneider} D. Schattschneider {\em M.C. Escher: Visions of Symmetry} second edition, Harry N. Abrams.

\bibitem{wang1961} H. Wang {\em Proving theorems by pattern 
recognition- II}, Bell System Technical Journal {\bf 40} (1961), 
1-42. 

\bibitem{washcrowe} D.K. Washburn and D.W. Crowe, {\em Symmetries of Culture: Theory and Practice of Plane Pattern Analysis}, University of Washington Press (1987).

\end{thebibliography}
\end{document}